\documentclass[a4paper,11pt]{amsart}
\usepackage[english]{babel}
\usepackage{amsmath,amssymb,amscd}
\usepackage{amsthm}
\usepackage{verbatim}
\usepackage[utf8]{inputenc}
\usepackage{amsaddr}
\usepackage{hyperref}
\usepackage{array}

\newtheorem{theorem}{Theorem}[section]

\numberwithin{equation}{section}

\begin{document}

		\title[rational non-integrability of the $N$-center problem]{A note on the rational non-integrability of the $N$-center problem  for almost all degrees of the singularities}
	
		\author{Eddaly Guerra-Velasco$^{1,2}$, Boris Percino-Figueroa$^2$, Russell-Aar\'on Qui\~nones-Estrella$^2$}
\email{edaly.velasco@unach.mx, eguerra@conahcyt.mx}
	
	\email{boris.percino@unach.mx}
	
	\email{rusell.quinones@unach.mx}
	\address{$^1$CONAHCyT, Av. de los Insurgentes Sur 1528, Benito Juarez, 03940, Ciudad de M\'exico, M\'exico.}	
	\address{$^2$Facultad de Ciencias en F\'isica y Matem\'aticas, Universidad Aut\'onoma de Chiapas, Carr. Emiliano Zapata km 8, Tuxtla Guti\'errez, 29050, Chiapas, M\'exico}
	
	\begin{abstract}
	In this article, we show that the $N$-center problem with rational weak and moderate forces is not rationally integrable for all but a finite number of values $\alpha\in(0,2)\cap \mathbb{Q}$, where $\alpha$ is the order of the singularities. We identify the remaining cases and provide the necessary conditions for their integrability. 
\end{abstract}


\maketitle

\section{Introduction}
In \cite{Sh}, M. Shibayama analyzes  the rational integrability of the Newtonian $N$-center problem; to be more specific, he considers the Hamiltonian system with three degrees of freedom, defined by the Hamiltonian function

\begin{equation}\label{ElH}
	H(q,p)=\frac12 \left\| p \right\|^2-\sum_{i=1}^{N}\frac{m_i}{\|q-c_i\|},
\end{equation}
where  $c_1,\dots,c_N\in \mathbb{R}^3$ are the fixed centers of the problem ($c_k=(c_k^1,c_k^2,c_k^3)$). By taking $l\in\mathbb{N}, l\leq N$ and a point $e=(e_1,e_2,e_3)\in\mathbb{C}^3$ such that
\begin{equation}
	(e_1-c^1_k)^2+(e_2-c^2_k)^2+(e_3-c_k^3)^2\begin{cases}
		=0, &k\leq l,\\
		\neq0, &k>l.
	\end{cases}, 
\end{equation}
he proves the following result (see \cite[Theorem 1]{Sh}).

\begin{theorem}

	If the $N$-center problem is given by the Hamiltonian function
	\eqref{ElH}
	is integrable, then for each $k\in\{1,2,3\}$, $\nu_k=\sqrt{\frac{25}{4}-\frac{4\rho_k}{C}}\in\frac12+\mathbb{Z}$, where $\rho_1,\rho_2,\rho_3$ are the eigenvalues of the matrix $A$ given by
	\begin{equation}\label{mat:B}
		A_{ij}=\sum_{k}\frac{3m_k(e_i-c_i^k)(e_j-c_j^k)}{[2q(t)\cdot(e-c_k)]^{\frac{3}{2}}}.
	\end{equation} 
	and $d$ satisfies the equation
	\[
	Cd=\sum_{k=1}^l\frac{m_k}{(2(e-c_k)\cdot d)^{\frac32}}(e-c_k),
	\]
	for some constant $C\in \mathbb{C}$.
\end{theorem}

\section{Rational weak and moderate forces}
In this manuscript, we notice that using the ideas, \emph{mutatis mutandis}, of Shibayama, this result can be generalized to a broader family of $N$-center problems to obtain a criterion on the integrability in terms of the degree of the singularities. To be more specific, consider the Hamiltonian system with $n$ degrees of freedom given by the Hamiltonian function
\begin{equation}\label{ElHaf}
	H_{\alpha}(q,p)=\frac12 \left\| p \right\|^2-\sum_{i=1}^{N}\frac{m_i}{\|q-c_i\|^{\alpha}},
\end{equation}
with centers  $c_1,\dots,c_N\in \mathbb{R}^n$ and $\alpha$ a rational number in $(0,2)$ called degree of the singularities.

In order to state the generalization, pick some $l\in\mathbb{N}, l\leq N$ and $e=(e_1,\dots,e_n)\in\mathbb{C}^n$ such that 
\[
(e_1-c_k^1)^2+\cdots+(e_n-c_k^n)^2=\begin{cases}
	=0, &k\leq l,\\
	\neq0, &k>l.
\end{cases}.
\]

Given a  vector $V=(v_1,\dots,v_n)\in\mathbb{C}^n$ define the constant
\begin{equation}\label{alg:rel}
	C:=C_V=\sum_{i=1}^{l}\frac{m_i}{\alpha[v_i(e-c_i]^{\frac{\alpha+2}{2}}}(e-c_i).
\end{equation}
Consider the system 
\begin{equation}\label{hamH0}
	\begin{split}
		\frac{dq}{dt}&=p\\
		\frac{dp}{dt}&=\sum_i \frac{m_i}{\frac{\alpha}{2}[2q\cdot(e-c_i)]^{\frac{\alpha+2}{2}}}(e-c_i),
	\end{split}
\end{equation}
it is straightforward that if $g(t)$ is a solution of 
\begin{equation}\label{kepler}
	\frac{d^2g(t)}{dt^2}=\frac{C_V}{(g(t))^{\frac{\alpha+2}{2}}}, 
\end{equation}
then 
\[
q(t)=g(t)V
\]
is a solution for \eqref{hamH0}. For such $q(t)$ define the $n\times n$ matrix $A$ with entries
\begin{equation}%\label{mat:B}
	A_{ij}=\sum_{k}\frac{(\alpha+2)m_k(e_i-c_i^k)(e_j-c_j^k)}{[2q(t)\cdot(e-c_k)]^{\frac{\alpha+2}{2}}}.
\end{equation} 
Finally,  consider the eigenvalues $a_k, k=1,\dots,n$ of the matrix $A$

As in Shimura's case, a key point of the proof is to compare the quantities 
\begin{equation}\label{relev:const}
	\begin{split}
		\lambda&=-\alpha-1\\
		\mu&=1-\frac12\alpha\\
		\nu_k&=\pm\sqrt{\Big(\frac32\alpha+1\Big)^2-4\frac{\alpha a_k}{C}}
	\end{split}
\end{equation} 
with the Schwarz's table in \cite{Ki} to have the following result:

\begin{theorem}[Main Theorem] \label{int:crit}
	If the $N$-center problem defined by $H_{\alpha}$ given in \eqref{ElHaf} is rationally integrable, then for every $ k\in\{1,\dots,n\}$, one of the cases in the following table holds for $\alpha$ and $a_k/C$.
	{\renewcommand{\arraystretch}{1.9}
		\begin{table}[h!]
			\centering
			\begin{tabular}{|c|c|c|}
				\hline  &$\alpha$&$a_k/C$\\
				\hline 1&$1$&$\sqrt{\frac{25}{4}-\frac{4 a_k}{C}}\in\frac12+\mathbb{Z}$  \\
				\hline 2&$\frac43$&$\sqrt{9-\frac{16 a_k}{3C}}\in(2\mathbb{Z}+1)\cup(\frac43+2\mathbb{Z})\cup\left(\frac65+2\mathbb{Z}\right)$\\
				\hline 3&$\frac23$& $\sqrt{1-\frac{2 a_k}{3C}}\in\frac12\mathbb{Z}\cup\left(\frac13+2\mathbb{Z}\right)\cup\left(\frac15+2\mathbb{Z}\right)$ \\
				\hline 4&$\frac12$&$\sqrt{\frac{49}{16}-\frac{2 a_k}{C}}\in\frac12+\mathbb{Z}$  \\
				\hline 5&$\frac32$&$\sqrt{\frac{169}{16}-\frac{6 a_k}{C}}\in\left(\frac13+\mathbb{Z}\right)\cup\left(\frac12+\mathbb{Z}\right)$\\
				\hline 6&$\frac45$&$\frac15\sqrt{121-\frac{80 a_k}{C}}\in\frac13+2\mathbb{Z}$   \\
				\hline 7&$\frac85$&$\sqrt{\frac{289}{25}-\frac{32 a_k}{5C}}\in\frac32+2\mathbb{Z}$  \\
				\hline
			\end{tabular}
			\caption{Integrability Criteria}
			\label{int:table}
	\end{table} }
	
	In particular, if $\alpha$ is not in the previous table, the N-center problem \eqref{ElHaf} is not integrable.
\end{theorem}

We finish by noticing that problem \eqref{ElHaf} has been studied by Bolotin \cite{Bo} and Koslov \cite{Ko}, who proved, by geometric methods, that the planar (\emph{i.e.} $n=2$) $N-$center problem, $N\geq 3$, is not integrable for moderate forces. By moderate forces, we mean that $1<\alpha<2$; recall that by weak forces, it is meant $0<\alpha<1$. This note extend this results to weak and moderate forces $\alpha\in ]0,2[\cap \mathbb{Q}$, for all but a finite number of $\alpha$ and arbitrary dimension $n$.

\end{document}